\documentclass[12pt,oneside]{amsart}%
\usepackage{amssymb,amsfonts}
\usepackage{amsfonts,amsmath,amsxtra,amsthm,amssymb,latexsym,mathrsfs,srcltx,dsfont,euscript}
\pdfoutput=1
\usepackage{amscd}

\usepackage[ruled,vlined,norelsize]{algorithm2e}
\usepackage{algorithmic}
\providecommand{\SetAlgoLined}{\SetLine}

\usepackage{graphicx} 
\usepackage{epstopdf}
\usepackage{caption}
\usepackage{subcaption}

\usepackage{hyperref} 

\usepackage[cp1251]{inputenc}
\usepackage[ukrainian, english]{babel}
\usepackage[left=2.5cm,right=2cm,
    top=2cm,bottom=2cm,bindingoffset=0cm]{geometry}

\newtheorem{theorem}{Theorem}

\numberwithin{equation}{section}
\numberwithin{theorem}{section}
\numberwithin{corollary}{section}
\numberwithin{remark}{section}
\numberwithin{lemma}{section}

\begin{document}

\address{Institute of  Mathematics of the NAS of Ukraine \\
3 Tereshchenkivs'ka St., Kyiv-4, 01601, Ukraine}
\author{Volodymyr L. Makarov}
\email[V.~L. Makarov]{makarov@imath.kiev.ua}
\author{Denys V. Dragunov}
\email[D.~V. Dragunov]{dragunovdenis@gmail.com}
\author{Danyil V. Bohdan}
\email[D.~V. Bohdan]{danyil.bohdan@gmail.com}

\title[Exponentially convergent numerical-analytical method\ldots]{Exponentially convergent numerical-analytical method for solving eigenvalue problems for singular differential operators}

\selectlanguage{english}

\begin{abstract}
The article develops and proves an exponentially convergent numerical-analytical method (the FD-method) for solving Sturm-Liouville problems with a singular Legendre operator and a singular potential. Obtained within are sufficient conditions for convergence of the method and a priori estimates of its accuracy.    A detailed algorithm for programmatic implementation of the FD-method is presented and compared with known algorithms (SLEIGN2).
\end{abstract}
\maketitle

\section{Introduction}
The results presented in this article constitute a logical extension and a generalization of the results in \cite{Klimenko_Ya} and \cite{Dragunov_Math_Comp}, which consider the subject of solving the following Sturm-Liouville problem:

\begin{equation}\label{eq_1}
    -\frac{d}{d x}\left[(1-x^{2})\frac{d u(x)}{d x}\right]+q(x)u(x)=\lambda u(x),\quad x\in(-1, 1),
\end{equation}
\begin{equation}\label{eq_1_1}
    \lim\limits_{x\rightarrow \pm 1}(1-x^{2})\frac{d u(x)}{d x}=0.
\end{equation}

Problems of this kind arise in applications when solving partial differential equations in spherical coordinates using  separation of variables, as is done, e.g., with \textit{hydrogen-molecule ion's equation} in \cite{Slavyanov} (see \cite[p. 167--170]{Slavyanov}).

To recall, articles \cite{Klimenko_Ya},  \cite{Dragunov_Math_Comp} develop and prove an exponentially convergent algorithm (an FD-method) for solving problem \eqref{eq_1}, \eqref{eq_1_1} for the case when the function $q(x)$ is of the class $Q^{0}[-1,1]$ of piecewise continuous functions that are bounded on the  closed interval $[-1,1]$ and have no more than a finite number of jump discontinuities. However,  \cite{Dragunov_Math_Comp} shows the results of applying the FD-method to problem \eqref{eq_1}, \eqref{eq_1_1} with the {\it potential} $q(x)=|x+1/3|^{1/2}+\ln(|x-1/3|),$ which clearly does not belong to the class $Q^{0}[-1,1].$ Despite the FD-method's convergence having not been proved for such problems the method turns out convergent. This fact has suggested to the authors of \cite{Dragunov_Math_Comp} that the sufficient convergence conditions for the FD-method for problems of type \eqref{eq_1}, \eqref{eq_1_1} can be weakened substantially, especially where they concern the smoothness of the function $q(x).$

The subject of this article is the Sturm-Liouville problem \eqref{eq_1}, \eqref{eq_1_1} with a function $q(x)$ from the space $L_{1, \rho}(-1, 1),$ $\rho=1/\sqrt{1-x^{2}},$ which contains functions $f(x),$ defined almost everywhere on the interval $(1-, 1)$ for which it holds that
\begin{equation}\label{eq_1_2}
    \|f\|_{1, \rho}=\int\limits_{-1}^{1}\frac{|f(x)|}{\sqrt{1-x^{2}}}d x<\infty.
\end{equation}
 Thus stated, the problem is a generalization of those considered in  \cite{Klimenko_Ya} and \cite{Dragunov_Math_Comp}, sufficiently so we could not apply the proof techniques used therein. Instead, to obtain the sufficient conditions for convergence a new approach was used based on an inequality for Legendre functions proposed by V. L. Makarov. This inequality (see Theorem \ref{Legendre_func_ineq_theorem}) follows from Theorem \ref{Elegant_theorem}, analogues of which the authors were unable to find. For this reason the aforementioned  theorems  are presented here with a detailed proof as novel and original results.

The article has the following structure: we start out by giving an outline of the FD-method in section 2 and applying it to the problem at hand. We proceed to prove a general auxiliary result in section 3. In section 4 we give a theoretical justification of the method as applied to the case at hand and obtain a proof of its convergence. We discuss the programmatic side of the question in section 5. Finally, we draw some conclusions about what has been done.

\section{The FD-method: algorithm}

We are going to construct a solving algorithm for problem \eqref{eq_1}, \eqref{eq_1_1} based on the general idea of the FD-method (see \cite{Makarov_Rossokhata_Review}).

It is easy to see that the differential operator $L[\cdot]$ defined by the equality
\begin{equation}\label{L_operator}
  L[u(x)]= \frac{d}{d x}\left[(1-x^{2})\frac{d u(x)}{d x}\right]-q(x)u(x)
\end{equation}
is self-adjoint  in the Hilbert space

\begin{equation}\label{hilbert_space}
  W=\left\{f(x)\in C^{2}(-1,1)\cap L_{2}(-1,1) \;|\; \lim\limits_{x\rightarrow \pm 1}(1-x^{2})f(x)=0;\; q(x)f(x)\in L_{2}(-1,1) \right\}
\end{equation}
 equipped with the common inner product
\begin{equation}\label{inner_product}
   <f,g>=\int\limits_{-1}^{1}f(x)g(x)dx
\end{equation}
 (see \cite[p. 55]{Al-Gwaiz_SL_theory}).
 This fact implies that there exists an increasing sequence of eigenvalues $\lambda_{0}<\lambda_{1}<\ldots < \lambda_{n}<\ldots$ and corresponding orthogonal eigenfunctions $u_{0}(x), u_{1}(x),\ldots, u_{n}(x), \ldots$ that satisfy equation \eqref{eq_1} and condition \eqref{eq_1_1}.

We are looking for the eigensolution $u_{n}(x), \lambda_{n}$ to eigenvalue problem \eqref{eq_1},\eqref{eq_1_1} in the form of a series
\begin{equation}\label{A_eq_1}
     u_{n}(x)=\sum\limits_{j=0}^{\infty}u^{(j)}_{n}(x),\quad \lambda_{n}=\sum\limits_{j=0}^{\infty}\lambda_{n}^{(j)},
\end{equation}
where the pair $u_{n}^{(j)}(x), \lambda_{n}^{(j)}$ can be found as the solution to the following system of recurrence problems:
\begin{equation}\label{A_eq_2}
    \frac{d}{d x}\left[(1-x^{2})\frac{d u_{n}^{(j)}(x)}{d x}\right]+\lambda_{n}^{(0)}u_{n}^{(j)}(x)=-\sum\limits_{i=0}^{j-1}\lambda_{n}^{(j-i)}u_{n}^{(i)}(x)+q(x)u_{n}^{(j-1)}(x),\; u_{n}^{(-1)}(x)\equiv 0,
\end{equation}
\begin{equation}\label{A_eq_2_1}
   \lim\limits_{x\rightarrow \pm 1}(1-x^{2})\frac{d u_{n}^{(j)}(x)}{d x}=0, \quad  j=0,1,2,\ldots.
\end{equation}
If we put $j=0$ in \eqref{A_eq_2} we obtain the equation for the {\it  basic problem}
\begin{equation}\label{A_eq_3}
   \frac{d}{d x}\left[(1-x^{2})\frac{d u_{n}^{(0)}(x)}{d x}\right]+\lambda_{n}^{(0)}u_{n}^{(0)}(x)=0.
\end{equation}
Taking into account that the eigenfunctions of operator $L[\cdot]$ \eqref{L_operator} are determined up to a multiplicative constant we impose an additional requirement on the solutions of the basic problem \eqref{A_eq_3}, \eqref{A_eq_2_1}:
\begin{equation}\label{norm_requirement}
    \int\limits_{-1}^{1}\left(u_{n}^{(0)}(x)\right)^{2}d x=1, \quad n=0, 1, 2,\ldots
\end{equation}

It is well known (see \cite[p. 121]{Bateman_HTF1}, \cite[p. 33]{King_Billingham_Diff_eqs}) that every solution $u_{n}^{(0)}(x)$ to equation \eqref{A_eq_3} (when $\lambda_{n}^{(0)}$ is fixed) can be represented through the Legendre functions $P_{\nu}(x),$ $Q_{\nu}(x):$
\begin{equation}\label{A_eq_4}
  u_{n}^{(0)}(x)=A P_{\nu}(x)+B Q_{\nu}(x), \quad A,B\in \mathbb{C},
\end{equation}
where $\nu$ is the solution of the algebraic equation $\nu(\nu+1)=\lambda_{n}^{(0)},$ i.e,
\begin{equation}\label{A_eq_5}
    \nu=-\frac{1}{2}\left(1\pm\sqrt{1+4\lambda_{n}^{(0)}}\right).
\end{equation}
Now taking into account the formulas that describe the behaviour of Legendre functions near the singular points $\pm 1$ (see \cite[p. 163--164]{Bateman_HTF1}) and the formulas that connect derivatives of the Legendre functions with the associated Legendre functions (see \cite[p. 148]{Bateman_HTF1}) we can easily compute that
\begin{eqnarray}\label{limits_eq}
& \lim\limits_{x\rightarrow -1}\left(1-x^{2}\right)\cfrac{d P_{\nu}(x)}{d x}=\cfrac{2\sin\left(\pi \nu\right)}{\pi}, \quad \lim\limits_{x\rightarrow 1}\left(1-x^{2}\right)\cfrac{d P_{\nu}(x)}{d x}=0, \nonumber\\
\\
& \lim\limits_{x\rightarrow -1}\left(1-x^{2}\right)\cfrac{d Q_{\nu}(x)}{d x}=\cos\left(\pi \nu\right), \quad \lim\limits_{x\rightarrow 1}\left(1-x^{2}\right)\cfrac{d Q_{\nu}(x)}{d x}=1.\nonumber
\end{eqnarray}
From equalities \eqref{limits_eq} it follows that function $u_{n}^{(0)}(x)$ \eqref{A_eq_4} satisfies condition \eqref{A_eq_2_1} if and only if $B=0$ and $\nu=n\in \mathbb{N}\cup\{0\},$ whereas condition \eqref{norm_requirement} leads us to the equality (see \cite[p. 42]{King_Billingham_Diff_eqs}) \begin{equation}\label{A_inv}
A^{-1}=\sqrt{\textstyle \int\limits_{-1}^{1}\left(P_{n}(x)\right)^{2}d x}=\sqrt{\frac{2}{2n+1}}.
\end{equation}
In the other words we have that the pairs
\begin{equation}\label{A_eq_6}
   u_{n}^{(0)}(x)=\sqrt{\frac{2n+1}{2}}P_{n}(x), \quad \lambda_{n}^{(0)}=n(n+1), \quad n=0,1,2,\ldots.
\end{equation}
represent eigensolutions of eigenvalue problem \eqref{A_eq_2_1}, \eqref{A_eq_3} and \eqref{norm_requirement}. Problems \eqref{A_eq_2}, \eqref{A_eq_2_1} for  $j=1,2,\ldots$ are solvable if and only if the functions
\begin{equation}\label{right_hand_side}
    F^{(j)}_{n}(x)=-\sum\limits_{i=0}^{j-1}\lambda_{n}^{(j-i)}u_{n}^{(i)}(x)+q(x)u_{n}^{(j-1)}(x), \;j=1,2,\ldots
\end{equation}
are orthogonal (in the sense of inner product \eqref{inner_product}) to the kernel space of the linear operator
$$L^{(0)}_{n}[u(x)]= \frac{d}{d x}\left[(1-x^{2})\frac{d u(x)}{d x}\right]+\lambda_{n}^{(0)}u(x),$$ i.e, to the function $u_{n}^{(0)}(x).$  This fact gives us a simple formula for finding $\lambda_{n}^{(j)},$ $j=1,2,\ldots$:
\begin{equation}\label{formula_for_lambda_j}
  \lambda_{n}^{(j)}=\int\limits_{-1}^{1}u_{n}^{(0)}(x)\left\{-\sum\limits_{i=1}^{j-1}\lambda_{n}^{(j-i)}u_{n}^{(i)}(x)+q(x)u_{n}^{(j-1)}(x)\right\} d x,
\end{equation}
whereas functions $u_{n}^{(j)}(x),$ $j=1,2,\ldots$ can be found via the variation of parameters formula (see, e.g., \cite[p. 8, 34]{King_Billingham_Diff_eqs})
\begin{equation}\label{formula_for_u_j}
  u_{n}^{(j)}(x)=c_{n}^{(j)}u_{n}^{(0)}(x)+\int\limits_{-1}^{x}K_{n}(x,\xi)F^{(j)}_{n}(\xi) d \xi,
\end{equation}
where
\begin{equation}\label{Cauchy_kernel}
   K_{n}(x, \xi)=P_{n}(x)Q_{n}(\xi)-Q_{n}(x)P_{n}(\xi)
\end{equation}
and constant $c_{n}^{(j)}\in \mathbb{R}$ can be chosen arbitrary. In a later section we will choose it to satisfy the orthogonality condition $\left<u_{n}^{(0)}(x), u_{n}^{(j)}(x)\right>=0.$

\section{Auxiliary results}
In what follows we will need the result stated below in the form of a theorem, which we consider to be quite elegant.
\begin{theorem}\label{Elegant_theorem}
   Suppose that $u_{I}(\theta)$ and $u_{II}(\theta)$ are a pair of solutions to the differential equation
   \begin{equation}\label{eq_3}
    \frac{d^{2} u(\theta)}{d\theta^{2}}+\phi(\theta)u(\theta)=0, \quad \theta\in \left(a, b\right),
   \end{equation}
   $$\phi(\theta)\in C^{1}(a,b), \phi(\theta)> 0, \forall \theta\in (a, b)$$
   that satisfy the following condition:
   \begin{equation}\label{eq_4}
    W(\theta)=u_{I}(\theta)u_{II}^{\prime}(\theta)-u^{\prime}_{I}(\theta)u_{II}(\theta)=1, \quad \forall \theta\in (a,b).
   \end{equation}
   If there exists a point $c\in (a, b)$ such that $\phi^{\prime}(\theta)\leq 0$ $\forall \theta\in (a,c]$ and $\phi^{\prime}(\theta)\geq 0$ $\forall \theta\in [c,b)$ then
\begin{equation}\label{eq_5}
   \left|v(\theta, \tilde{\theta})\right| \leq \sqrt{2\phi^{-1}(c)}, \forall \theta, \tilde{\theta}\in (a, b),
   \end{equation}
   $$v(\theta, \tilde{\theta})  \stackrel{def}{=} u_{I}(\theta)u_{II}(\tilde{\theta})-u_{I}(\tilde{\theta})u_{II}(\theta).$$
   If $\phi^{\prime}(\theta)\leq 0$ or $\phi^{\prime}(\theta)\geq 0$  $\forall \theta\in (a,b)$ then
   \begin{equation}\label{eq_5_1}
   \left|v(\theta, \tilde{\theta})\right| \leq \max\left\{\sqrt{\phi^{-1}(\theta)}, \sqrt{\phi^{-1}(\tilde{\theta})}\right\}, \forall \theta, \tilde{\theta}\in (a, b).
   \end{equation}
\end{theorem}

Before proceeding to the proof of Theorem \ref{Elegant_theorem} we should emphasize that the main idea of the theorem was evoked by the Theorem of Sonin (see \cite[p. 166]{Szego_orto_poly}).

\begin{proof}
Suppose that the conditions of Theorem \ref{Elegant_theorem} are fulfilled and for some $c\in (a,b)$ we have that
\begin{equation}\label{c_point_inequalities}
\phi^{\prime}(\theta)\leq 0\;\forall \theta\in (a,c],\;\;\phi^{\prime}(\theta)\geq 0\;\forall \theta\in [c,b).
\end{equation}
In such a case the auxiliary function
\begin{equation}\label{eq_6}
    f_{1}(\theta, \tilde{\theta})=v^{2}(\theta, \tilde{\theta})+\phi^{-1}(\theta)\left(\frac{\partial v(\theta, \tilde{\theta})}{\partial \theta}\right)^{2}
\end{equation}
is non-decreasing on $(a,c]$ and non-increasing on $[c, b)$ with respect to its argument $\theta,$ i.e.,
\begin{equation}\label{eq_6_iii}
\frac{\partial f_{1}(\theta, \tilde{\theta})}{\partial \theta}\geq 0, \; \forall \theta\in (a, c];\quad \frac{\partial f_{1}(\theta, \tilde{\theta})}{\partial \theta}\leq 0,\; \forall \theta\in [c, b),\; \forall \tilde{\theta}\in(a,b).
\end{equation}
The latter fact easily follows from the equality
$$\frac{\partial f_{1}(\theta, \tilde{\theta})}{\partial \theta}=2\frac{\partial v(\theta, \tilde{\theta})}{\partial \theta}\left(v(\theta, \tilde{\theta})+\phi^{-1}(\theta)\frac{\partial^{2} v(\theta, \tilde{\theta})}{\partial \theta^{2}}-\frac{1}{2}\frac{\phi^{\prime}(\theta)}{\phi^{2}(\theta)}\left(\frac{\partial v(\theta, \tilde{\theta})}{\partial \theta}\right)\right)=$$
\begin{equation}\label{eq_6_i}
=-\frac{\phi^{\prime}(\theta)}{\phi^{2}(\theta)}\left(\frac{\partial v(\theta, \tilde{\theta})}{\partial \theta}\right)^{2},\quad \forall \theta,\tilde{\theta}\in (a,b)
\end{equation}
and inequalities \eqref{c_point_inequalities}.
In much the same way it is easy to verify that
\begin{equation}\label{eq_6_ii}
\frac{\partial f_{2}(\theta, \tilde{\theta})}{\partial \tilde{\theta}}\geq 0, \; \forall \tilde{\theta}\in (a, c];\quad \frac{\partial f_{2}(\theta, \tilde{\theta})}{\partial \tilde{\theta}}\leq 0,\; \forall \tilde{\theta}\in [c, b),\; \forall \theta\in(a,b).
\end{equation}
where
\begin{equation}\label{eq_7}
    f_{2}(\theta, \tilde{\theta})=v^{2}(\theta, \tilde{\theta})+\phi^{-1}(\tilde{\theta})\left(\frac{\partial v(\theta, \tilde{\theta})}{\partial \tilde{\theta}}\right)^{2}.
\end{equation}

\begin{figure}[t]
	\centering
	\includegraphics[width=0.45\linewidth]{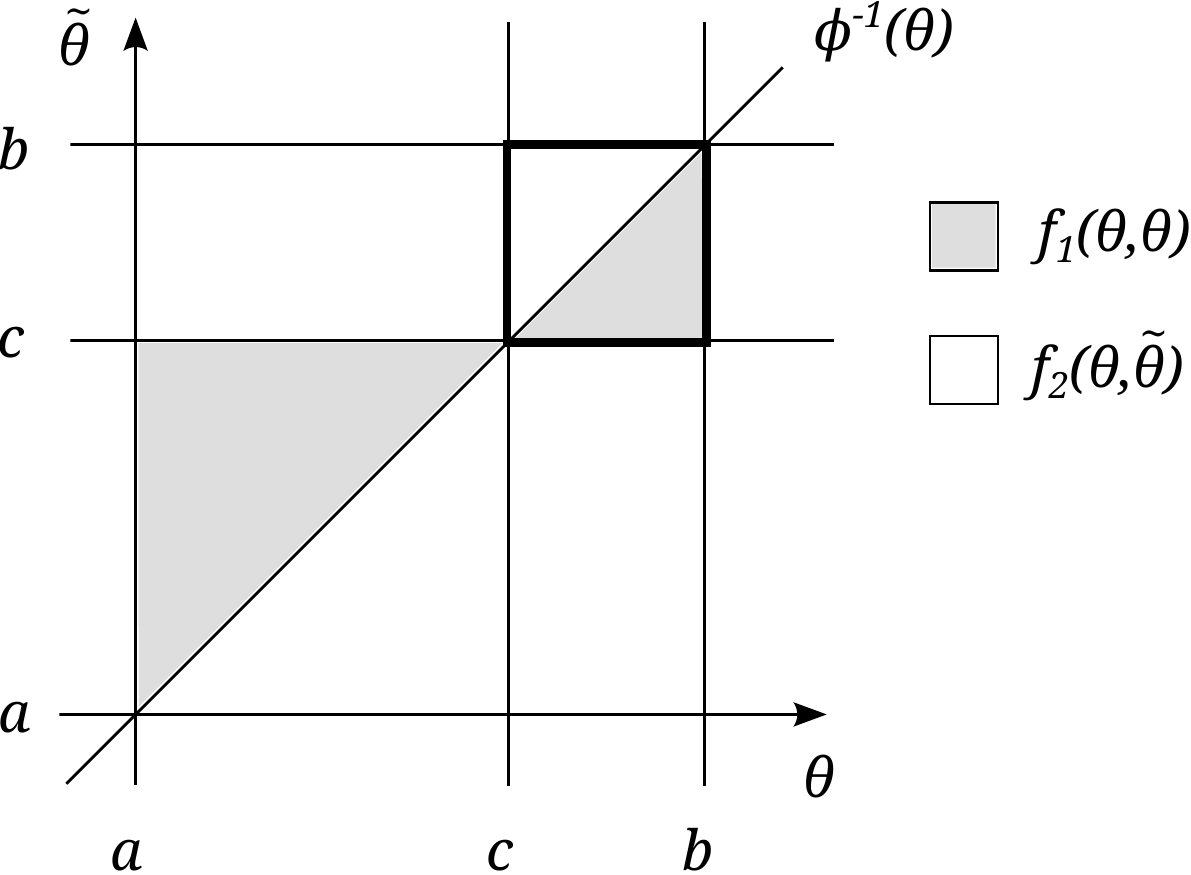}	
	\caption{A graph of $f(\theta, \tilde{\theta})$, $\bar{f}(\theta, \tilde{\theta})$ on $(a,b)\times(a,b)$. The domain of $\bar{f}(\theta, \tilde{\theta})$ is indicated with a bold-lined square. It may be of interest to note that on $(a,c)^2\cup(c,b)^2$ the right side of estimate (\ref{eq_5}) need not contain $\sqrt{2}$; see formulas (\ref{eq_12}), (\ref{eq_13}).
	}
	\label{th31fig}	
\end{figure}

Let us consider the function $f(\theta, \tilde{\theta})$ defined in the following way:
\begin{equation}\label{eq_8}
    f(\theta, \tilde{\theta})=\left\{
                \begin{array}{cc}
                  f_{1}(\theta, \tilde{\theta}) & when \quad \theta\leq \tilde{\theta}; \;  \theta,\tilde{\theta} \in (a, c],\\
                  f_{2}(\theta, \tilde{\theta}) & when \quad \theta\geq \tilde{\theta}; \;  \theta,\tilde{\theta} \in (a, c].\\
                \end{array}
              \right.
\end{equation}
From expressions \eqref{eq_6} and \eqref{eq_7} it follows that $f_{1}(\theta, \theta)=f_{2}(\theta, \theta)=\phi^{-1}(\theta).$ The latter fact means that the function $f(\theta, \tilde{\theta})$ \eqref{eq_8} is well defined and continuous on $(a, c]^{2}.$
Expressions \eqref{eq_6} and \eqref{eq_7} together with inequalities \eqref{eq_6_iii}, \eqref{eq_6_ii} and identity \eqref{eq_4} lead us to the inequalities
\begin{equation}\label{eq_10}
    f_{1}(\theta, \tilde{\theta})\leq f_{1}(\tilde{\theta}, \tilde{\theta})=\phi^{-1}(\tilde{\theta}), \quad \forall \theta,\tilde{\theta} \in (a, c], \; \theta\leq \tilde{\theta},
\end{equation}
\begin{equation}\label{eq_11}
    f_{2}(\theta, \tilde{\theta})\leq f_{2}(\theta, \theta)=\phi^{-1}(\theta), \quad \forall \theta,\tilde{\theta} \in (a, c], \; \theta\geq \tilde{\theta}.
\end{equation}
Taking into account expressions \eqref{eq_6}, \eqref{eq_7} from inequalities \eqref{eq_10}, \eqref{eq_11} we can deduce that
\begin{equation}\label{eq_12}
v^{2}(\theta, \tilde{\theta})\leq f(\theta, \tilde{\theta})\leq\phi^{-1}(c), \quad \forall \theta,\tilde{\theta} \in (a, c].
\end{equation}
Through applying nearly identical reasoning to the function
\begin{equation*}
    \bar{f}(\theta, \tilde{\theta})=\left\{
                \begin{array}{cc}
                  f_{2}(\theta, \tilde{\theta}) & when \quad \theta\leq \tilde{\theta}, \;  \theta,\tilde{\theta}\in [c, b),\\
                  f_{1}(\theta, \tilde{\theta}) & when \quad \theta\geq \tilde{\theta}, \;  \theta,\tilde{\theta}\in [c, b).\\
                \end{array}
              \right.
\end{equation*}
we can get the inequalities
\begin{equation}\label{eq_13}
v^{2}(\theta, \tilde{\theta})\leq \bar{f}(\theta, \tilde{\theta})\leq \phi^{-1}(c), \quad \forall \theta,\tilde{\theta} \in [c, b).
\end{equation}

Now let us return to the function $f_{1}(\theta, \tilde{\theta})$ \eqref{eq_6}. From inequalities \eqref{eq_6_iii} it follows that
\begin{equation}\label{eq_14}
    v^{2}(\theta,\tilde{\theta})\leq f_{1}(\theta, \tilde{\theta})\leq f_{1}(c, \tilde{\theta})=v^{2}(c,\tilde{\theta})+\phi^{-1}(c)w^{2}(\tilde{\theta}), \quad \theta,\tilde{\theta}\in (a,b)
\end{equation}
where $w(\tilde{\theta})=u^{\prime}_{I}(c)u_{II}(\tilde{\theta})-u_{I}(\tilde{\theta})u^{\prime}_{II}(c).$
Taking into account inequalities \eqref{eq_12} and \eqref{eq_13} we can proceed estimating $ v^{2}(\theta,\tilde{\theta})$ as follows:
\begin{equation}\label{eq_15}
  v^{2}(\theta,\tilde{\theta})\leq \phi^{-1}(c)(1+w^{2}(\tilde{\theta})),\quad \theta,\tilde{\theta}\in(a,b).
\end{equation}
It is not hard to verify that
$$w^{2}(\tilde{\theta})\leq w^{2}(\tilde{\theta})+\phi^{-1}(\tilde{\theta})\left(\frac{d w(\tilde{\theta})}{d\tilde{\theta}}\right)^{2}\leq w^{2}(c)+\phi^{-1}(c)\left.\left(\frac{d w(\tilde{\theta})}{d\tilde{\theta}}\right)^{2}\right|_{\tilde{\theta}=c}=1.$$
Combining the latter inequality with inequality \eqref{eq_15} we arrive at sought inequality \eqref{eq_5}.

Inequality \eqref{eq_5_1} can be easily derived from inequalities \eqref{eq_10}, \eqref{eq_11} as a limit case when $c\rightarrow b.$

The proof is complete.
\end{proof}

Using Theorem \ref{Elegant_theorem} we can obtain a curious and useful inequality pertaining to the Legendre functions.

It is well known that the Legendre functions $P_{\nu}(x)$ and $Q_{\nu}(x)$ are two linearly independent solutions to the Legendre differential equation (see \cite[p. 121]{Bateman_HTF1}):
\begin{equation}\label{eq_16}
    \frac{d}{d x}\left[(1-x^{2})\frac{d y(x)}{d x}\right]+\nu(\nu+1)y(x)=0,\quad x\in (-1, 1).
\end{equation}
Furthermore, the functions $P_{\nu}(x)$ and $Q_{\nu}(x)$ possess the property
\begin{equation}\label{eq_17}
    (1-x^{2})(P^{\prime}_{\nu}(x)Q_{\nu}(x)-P_{\nu}(x)Q^{\prime}_{\nu}(x))= 1, \quad \forall x\in (-1, 1).
\end{equation}
It is also known (see \cite[p. 67]{Szego_orto_poly}) that equation \eqref{eq_16} can be rewritten in the equivalent form \eqref{eq_3} with
\begin{equation}\label{eq_18}
   \phi(\theta)=(2\sin(\theta))^{-2}+(\nu+{1}/{2})^{2},
\end{equation}
$$u(\theta)=\sqrt{\sin(\theta)}y(\cos(\theta)), \quad a=0,\; b=\pi.$$
In the other words, we have that functions
\begin{equation}\label{eq_19}
    u_{I}(\theta)=\sqrt{{\sin(\theta)}}P_{\nu}(\cos(\theta)),\quad u_{II}(\theta)=\sqrt{{\sin(\theta)}}Q_{\nu}(\cos(\theta))
\end{equation}
satisfy equations \eqref{eq_3}, \eqref{eq_18} and identity \eqref{eq_4}, which is equivalent to identity  \eqref{eq_17}. Also, it is easy to see that the function $\phi(\theta)$ \eqref{eq_18} fulfils all the requirements of Theorem \ref{Elegant_theorem} with $c=\pi/2.$ Therefore, Theorem \ref{Elegant_theorem} provides us with the estimation
$$\sqrt{{\sin(\theta)}{\sin(\tilde{\theta})}}\cdot\left|P_{\nu}(\cos(\theta))Q_{\nu}(\cos(\tilde{\theta})-P_{\nu}(\cos(\tilde{\theta}))Q_{\nu}(\cos({\theta})\right|\leq\sqrt{2\phi^{-1}(\pi/2)}\leq \sqrt{\frac{2}{\frac{1}{4}+\left(\nu+\frac{1}{2}\right)^{2}}}, $$
$\forall \theta, \tilde{\theta}\in (0, \pi)$ and the following corollary:
\begin{theorem}\label{Legendre_func_ineq_theorem}
  For every $\nu\in \mathbb{R}$ the inequality
  \begin{equation}\label{eq_20}
    \sqrt[4]{(1-x^2)(1-\xi^2)}\left|P_{\nu}(x)Q_{\nu}(\xi)-P_{\nu}(\xi)Q_{\nu}(x)\right|\leq \sqrt{\frac{2}{\frac{1}{4}+\left(\nu+\frac{1}{2}\right)^{2}}}
  \end{equation}
  holds true for all $x,\xi\in (-1, 1).$
\end{theorem}

\section{The FD-method: theoretical justification}

In this section we are going to investigate the question of convergence of the proposed FD-method, i.e, to find the sufficient conditions that provide the convergence of series \eqref{A_eq_1}.

Let us consider a general eigenvalue problem

\begin{equation}\label{TJ_eq_1}
    \frac{d}{d x}\left[(1-x^{2})\frac{d u(x,\tau)}{d x}\right]-\tau q(x) u(x,\tau)=-\lambda(\tau)u(x,\tau), \; x\in (-1,1),\; \tau\in [0,1],
\end{equation}
\begin{equation}\label{TJ_eq_2}
    \lim\limits_{x\rightarrow \pm 1}(1-x^{2})\frac{d u(x,\tau)}{d x}=0,\; \forall \tau\in [0,1].
\end{equation}
The problem (\ref{eq_1}), (\ref{eq_1_1}) is its partial case for $\tau=1$. If we suppose that the eigenvalue $\lambda_{n}(\tau)$ and the corresponding eigenfunction $u_{n}(x,\tau)$ can be expressed in the form of a series
$$\lambda_{n}(\tau)=\sum\limits_{i=0}^{\infty}\lambda_{n}^{(i)}\tau^{i},\quad u_{n}(x,\tau)=\sum\limits_{i=0}^{\infty}u^{(i)}(x)\tau^{i}\quad \forall x\in (-1,1), \; \tau\in [0,1]$$
and the differential formulas
$$\frac{d u_{n}(x,\tau)}{d x}=\sum\limits_{i=0}^{\infty}\frac{d u_{n}^{(i)}(x)}{d x}\tau^{i},\quad \frac{d^{2} u_{n}(x,\tau)}{d x^{2}}=\sum\limits_{i=0}^{\infty}\frac{d^{2} u_{n}^{(i)}(x)}{d x^{2}}\tau^{i}\quad \forall x\in (-1,1), \; \tau\in [0,1]$$
hold we immediately arrive at the conclusion that the unknown coefficients $\lambda_{n}^{(i)}, u_{n}^{(i)}(x), \; i=0,1,2,\ldots$ can be found as solutions to problems \eqref{A_eq_2}, \eqref{A_eq_2_1}. To justify formulas \eqref{A_eq_1}  we only need to mention that if we set $\tau=1$ problem \eqref{TJ_eq_1}, \eqref{TJ_eq_2} will be reduced to problem \eqref{eq_1}, \eqref{eq_1_1}.

Now let us go back to formula \eqref{formula_for_u_j}. Without loss of generality we can obtain the values of $c_{n}^{(j)}$ using the orthogonality condition:
\begin{equation}\label{TJ_eq_3}
    c^{j}_{n}=-\int\limits_{-1}^{1}u^{(0)}_{n}(x)\int\limits_{-1}^{x}K_{n}(x,\xi)F^{(j)}_{n}(\xi)d\xi d x.
\end{equation}
It is not hard to verify that if $c_{n}^{(j)}$ is found according to formula \eqref{TJ_eq_3} then
$$\left<u_{n}^{(0)}(x), u_{n}^{(j)}(x)\right>\stackrel{def}{=}\int\limits_{-1}^{1}u_{n}^{(0)}u_{n}^{(j)}(x)d x=0, \quad \forall j\in \mathbb{N}$$
and formula \eqref{formula_for_lambda_j} can be substantially simplified:
\begin{equation}\label{TJ_eq_4}
    \lambda_{n}^{(j)}=\int\limits_{-1}^{1}q(x)u_{n}^{(0)}(x) u_{n}^{(j-1)}(x) d, \quad j\in \mathbb{N}.
\end{equation}

Using the norm $\|\cdot\|_{1,\rho}$ introduced in \eqref{eq_1_2} and formula \eqref{TJ_eq_4} we can estimate $|\lambda_{n}^{(j)}|$ as follows:
\begin{equation}\label{TJ_eq_5}
    \left|\lambda_{n}^{(j)}\right|=\left|\int\limits_{-1}^{1}\frac{q(x)}{\sqrt{1-x^{2}}}\sqrt[4]{1-x^{2}}u_{n}^{(j-1)}(x)\sqrt[4]{1-x^{2}}u_{n}^{(0)}(x)d x\right|\leq
\end{equation}
$$\leq \|q\|_{1, \rho} \|u_{n}^{(j-1)}\|_{\infty, 1/\sqrt{\rho}}\|u_{n}^{(0)}\|_{\infty, 1/\sqrt{\rho}},$$
where $\|f\|_{\infty, 1/\sqrt{\rho}}\stackrel{def}{=}\max\limits_{x\in [-1,1]}|f(x)\sqrt{\rho}|=\max\limits_{x\in [-1,1]}|\sqrt[4]{1-x^{2}}f(x)|.$ Theorem 7.3.3 from \cite{Szego_orto_poly} allows us to estimate $\|u_{n}^{(0)}\|_{\infty, 1/\sqrt{\rho}}$ as follows:
\begin{equation}\label{TJ_eq_6}
    \|u_{n}^{(0)}\|_{\infty, 1/\sqrt{\rho}}=\sqrt{n+1/2}\max\limits_{0\leq\theta\leq\pi}\sqrt{\sin(\theta)}\left|P_{n}(\cos(\theta))\right|\leq\sqrt{\frac{2(n+1/2)}{\pi n}}.
\end{equation}

Combining the latter inequality with \eqref{TJ_eq_5} we get the estimation
\begin{equation}\label{TJ_eq_7}
    \left|\lambda_{n}^{(j)}\right|\leq \sqrt{\frac{2(n+1/2)}{\pi n}}\|q\|_{1,\rho}\|u_{n}^{(j-1)}\|_{\infty, 1/\sqrt{\rho}}.
\end{equation}
Using estimation \eqref{TJ_eq_6} and formula \eqref{TJ_eq_3} we can estimate $|c_{n}^{(j)}|$ as follows:
$$\left|c_{n}^{(j)}\right|=\left|\int\limits_{-1}^{1}\frac{1}{\sqrt{1-x^{2}}}\sqrt[4]{1-x^{2}}u_{n}^{(0)}(x)\int\limits_{-1}^{x}\frac{\sqrt[4]{(1-x^{2})(1-\xi^{2})}}{\sqrt{1-\xi^{2}}}K_{n}(x,\xi)\right.\times$$
$$\times\left.\left[-\sum\limits_{i=0}^{j-1}\lambda_{n}^{(j-i)}\sqrt[4]{1-\xi^{2}}u_{n}^{(i)}(\xi)+q(\xi)\sqrt[4]{1-\xi^{2}}u_{n}^{j-1}(\xi)\right]d\xi d x\right|\leq$$
$$\leq \frac{\sqrt{2}\pi^{2}}{n+1/2}\sqrt{\frac{2(n+1/2)}{\pi n}}\left[\sum\limits_{i=0}^{j-1}|\lambda_{n}^{(j-i)}|\|u_{n}^{(i)}\|_{\infty, 1/\sqrt{\rho}}+\|q\|_{1,\rho}\|u_{n}^{(j-1)}\|_{\infty, 1/\sqrt{\rho}}\right] =$$
$$= \frac{2\pi\sqrt{\pi}}{\sqrt{n}\sqrt{n+1/2}}\left[\sum\limits_{i=0}^{j-1}|\lambda_{n}^{(j-i)}|\|u_{n}^{(i)}\|_{\infty, 1/\sqrt{\rho}}+\|q\|_{1,\rho}\|u_{n}^{(j-1)}\|_{\infty, 1/\sqrt{\rho}}\right].$$
To obtain the latter inequality we used the evident equality
$$\int\limits_{-1}^{1}\frac{d x}{\sqrt{1-x^{2}}}=\pi$$
and the result of Theorem \ref{Legendre_func_ineq_theorem}.

Now we are in a position to estimate $\|u_{n}^{(j)}\|_{\infty, 1/\sqrt{\rho}}.$ We can do this in the following way (see formula (\ref{formula_for_u_j})):
$$\|u_{n}^{(j)}\|_{\infty, 1/\sqrt{\rho}}\leq |c_{n}^{(j)}|\sqrt{\frac{2(n+1/2)}{\pi n}}+\max\limits_{x\in [-1,1]}\left\{\int\limits_{-1}^{x}\frac{\sqrt[4]{1-x^{2}}\sqrt[4]{1-\xi^{2}}}{\sqrt{1-\xi^{2}}}|K_{n}(x, \xi)|d\xi\times\right.$$
$$\left.\times\left[\sum\limits_{i=0}^{j-1}|\lambda_{n}^{(j-i)}|\|u_{n}^{(i)}\|_{\infty, 1/\sqrt{\rho}}+\|q\|_{1,\rho}\|u_{n}^{(j-1)}\|_{\infty, 1/\sqrt{\rho}}\right] \right\}\leq$$
$$\leq |c_{n}^{(j)}|\sqrt{\frac{2(n+1/2)}{\pi n}}+\frac{\sqrt{2}\pi}{n+1/2}\left[\sum\limits_{i=0}^{j-1}|\lambda_{n}^{(j-i)}|\|u_{n}^{(i)}\|_{\infty, 1/\sqrt{\rho}}+\|q\|_{1,\rho}\|u_{n}^{(j-1)}\|_{\infty, 1/\sqrt{\rho}}\right]=$$
\begin{equation}\label{TJ_eq_8}
    =\sqrt{2}\pi\left(\frac{3n+1}{n(n+1/2)}\right)\left[\sum\limits_{i=0}^{j-1}|\lambda_{n}^{(j-i)}|\|u_{n}^{(i)}\|_{\infty, 1/\sqrt{\rho}}+\|q\|_{1,\rho}\|u_{n}^{(j-1)}\|_{\infty, 1/\sqrt{\rho}}\right].
\end{equation}
Combining inequalities \eqref{TJ_eq_7} and \eqref{TJ_eq_8} we arrive at the following estimate:
\begin{equation}\label{TJ_eq_9}
    \|u_{n}^{(j)}\|_{\infty, 1/\sqrt{\rho}}\leq \alpha_{n}\left[\beta_{n}\sum\limits_{i=0}^{j-1}\|u_{n}^{(j-i-1)}\|_{\infty, 1/\sqrt{\rho}}\|u_{n}^{(i)}\|_{\infty, 1/\sqrt{\rho}}+\|u_{n}^{(j-1)}\|_{\infty, 1/\sqrt{\rho}}\right],
\end{equation}
where $$\alpha_{n}= \alpha_{n}(n) = \frac{\sqrt{2}\pi (3n+1)}{n(n+1/2)}\|q\|_{1,\rho}\leq\frac{3\sqrt{2}\pi}{n}\|q\|_{1,\rho},$$

$$\beta_{n}= \beta_{n}(n) =\sqrt{\frac{2(n+1/2)}{\pi n}}\leq\sqrt{\frac{3}{\pi}}<1.$$

 Using substitution
 \begin{equation}\label{TJ_eq_10}
    \|u_{n}^{j}\|_{\infty, 1/\sqrt{\rho}}=\alpha_{n}^{j}v_{j}
 \end{equation}
 we can rewrite inequality \eqref{TJ_eq_9} in the form of
 \begin{equation}\label{TJ_eq_11}
    v_{j}\leq \sum\limits_{i=0}^{j-1}v_{i}v_{j-i-1}+v_{j-1}, \quad j=1,2,\ldots,\; v_{0}=\|u_{n}^{(0)}\|_{\infty, 1/\sqrt{\rho}}.
 \end{equation}

 Let us consider a sequence of positive real numbers $\{V_{i}\}_{i=0,1,\ldots}$ defined by the recurrence formula
 \begin{equation}\label{TJ_eq_12}
    V_{j+1}=\sum\limits_{i=0}^{j}V_{i}V_{j-i}+V_{j}, \quad j=0,1,2,\ldots,\; V_{0}=1.
 \end{equation}

 Comparing \eqref{TJ_eq_11} with \eqref{TJ_eq_12} and taking into account inequality \eqref{TJ_eq_6} we can arrive at the conclusion that
 \begin{equation}\label{TJ_eq_13}
    v_{j}\leq V_{j},  j=0,1,2,\ldots.
 \end{equation}
 Recall that $\|u_{n}^{j}\|_{\infty, 1/\sqrt{\rho}}\leq \alpha_{n}^{j}V_{j}=\alpha_{n}^{j}V_{j}(n).$ If for some $n=n_{0}\in \mathbb{N}$ the series $\sum\limits_{i=0}^{\infty}\alpha_{n}^{j}V_{j}(n)$ is convergent then according the inequalities \eqref{TJ_eq_13}, \eqref{TJ_eq_5} and equality \eqref{TJ_eq_10} the series \eqref{A_eq_1} are convergent, i.e., the FD-method is convergent. Now we are going to find the smallest $n_{0}$ of the kind mentioned above. For this purpose let us consider the series
 \begin{equation}\label{TJ_eq_14}
    f(z)=\sum\limits_{j=0}^{\infty} z^{j}V_{j}
 \end{equation}
 and find its radius of convergence.

 Taking into account recurrence equalities \eqref{TJ_eq_12} one can verify that function $f(z)$ \eqref{TJ_eq_14} satisfies the functional equation
 $$f(z)=zf^{2}(z)+z f(z)+1$$ or, in a more convenient form,
\begin{equation}\label{TJ_eq_15}
 zf^{2}(z)+(z-1)f(x)+1=0.
\end{equation}
Solving equation \eqref{TJ_eq_15} with respect to function $f(z)$ we obtain
\begin{equation}\label{TJ_eq_16}
    f(z)=\frac{1}{2z}\left(1-z-\sqrt{1-\frac{z}{\gamma}}\sqrt{1-\gamma z}\right), \;\; \gamma=3 - 2\sqrt{2}.
\end{equation}
From formula \eqref{TJ_eq_16} we see that the radius of convergence $R$ for series \eqref{TJ_eq_14} is equivalent to $\gamma:$
\begin{equation}
    R=\gamma=3 - 2\sqrt{2}.
\end{equation}
Thus, if
\begin{equation}\label{TJ_eq_18}
   3 - 2\sqrt{2} \geq \frac{3\sqrt{2}\pi}{n}\|q\|_{1,\rho} \geq \alpha_{n}
\end{equation}
 the FD-method is convergent.

For a sufficiently large $n$ inequality \eqref{TJ_eq_18} will always be satisfied. This means that for sufficiently large values of $n$ the FD-method will always be convergent. To be specific, the FD-method will be convergent for all $n> n_{0},$ where  $\mathbb{N} \ni n_{0} \ge \frac{3\sqrt{2}\pi}{3 - 2\sqrt{2}}\|q\|_{1,\rho}$.

Furthermore, formula \eqref{TJ_eq_16} allows us to find the coefficients $V_{j},$ $j\in \mathbb{N}$ explicitly. For this purpose we need to expand the right-hand side of formula \eqref{TJ_eq_16} into a power series with respect to $z:$
$$f(z)=\frac{1}{2z}\left(1-z-\left[\sqrt{\gamma}-\frac{z}{2\sqrt{\gamma}}-\sum\limits_{p=2}^{\infty}\frac{(2p-3)!!}{(2p)!!}\gamma^{1/2-p}z^{p}\right]\times \right.$$
$$\times\left.\left[\frac{1}{\sqrt{\gamma}}-\frac{\sqrt{\gamma}z}{2}-\sum\limits_{p=2}^{\infty}\frac{(2p-3)!!}{(2p)!!}\gamma^{p-1/2}z^{p}\right]\right)=$$
$$=\frac{1}{2z}\left(-z+\sum\limits_{j=1}^{\infty}z^{j}\left[\frac{(2j-3)!!}{(2j)!!}(\gamma^{j}+\gamma^{-j})-\sum\limits_{p=1}^{j-1}\frac{(2p-3)!!}{(2p)!!}\frac{(2j-2p-3)!!}{(2j-2p)!!}\gamma^{2p-j}\right]\right)=$$
\begin{equation}\label{TJ_eq_20}
=1+\frac{1}{2}\sum\limits_{j=2}^{\infty}z^{j-1}\left[\frac{(2j-3)!!}{(2j)!!}(\gamma^{j}+\gamma^{-j})-\sum\limits_{p=1}^{j-1}\frac{(2p-3)!!}{(2p)!!}\frac{(2j-2p-3)!!}{(2j-2p)!!}\gamma^{2p-j}\right].
\end{equation}
Here we define $(2p)!!$ as $2\times 4\times\ldots \times 2p$ and $(2p+1)!!$ as $1\times 3\times\ldots \times (2p+1),$ $(-1)!!\stackrel{def}{=}1.$

From \eqref{TJ_eq_20} we have that
\begin{equation}\label{TJ_eq_21}
    V_{j-1}=\frac{1}{2}\left[\frac{(2j-3)!!}{(2j)!!}(\gamma^{j}+\gamma^{-j})-\sum\limits_{p=1}^{j-1}\frac{(2p-3)!!}{(2p)!!}\frac{(2j-2p-3)!!}{(2j-2p)!!}\gamma^{2p-j}\right], \; j=2,3,\ldots.
\end{equation}
Using the fact that $V_{j}\geq 0$ and Stirling's formula we can estimate $V_{j-1}$ in the following way:
$$V_{j-1}\leq \frac{(2j-3)!!}{2(2j)!!}=\frac{(2j-1)!!}{2(2j-1)(2j)!!}=\frac{(2j)!}{2(2j-1)((2j)!!)^{2}}=\frac{(2j)!}{2^{2j+1}(2j-1)(j!)^{2}}<$$
\begin{equation}\label{TJ_eq_22}
<\frac{2\sqrt{\pi j}(2j)^{2j}e^{-2j+1/(24 j)}}{2^{2j+1}(2j-1)(\sqrt{2\pi j}j^{j}e^{-j})^{2}}=\frac{e^{1/(24 j)}}{2(2j-1)\sqrt{\pi j}}<\frac{1}{(2j-1)\sqrt{\pi j}}.
\end{equation}
Using inequalities \eqref{TJ_eq_7}, \eqref{TJ_eq_13} and \eqref{TJ_eq_22} together with equality \eqref{TJ_eq_10} we can easily estimate $\|u_{n}^{(j)}\|_{\infty, 1/\sqrt{\rho}}$ and $|\lambda_{n}^{(j)}|:$
\begin{equation}\label{TJ_eq_23}
   \|u_{n}^{(j)}\|_{\infty, 1/\sqrt{\rho}}\leq \alpha_{n}^{j}V_{j}\leq \left(\frac{3\sqrt{2}\pi}{n}\|q\|_{1,\rho}\right)^j \frac{1}{(2j+1)\sqrt{\pi(j+1)}},
\end{equation}
\begin{equation}\label{TJ_eq_24}
   \left|\lambda_{n}^{(j)}\right|\leq \left(\frac{3\sqrt{2}\pi}{n}\right)^{j-1}(\|q\|_{1,\rho})^{j} \frac{1}{(2j-1)\sqrt{\pi j}}
\end{equation}

Inequalities \eqref{TJ_eq_23}, \eqref{TJ_eq_24}
now allow us to formulate the theorem about convergence of the FD-method.

\begin{theorem}
\label{conv_theorem}
 Let $$n_{0} = \left[\frac{3\sqrt{2}\pi}{3 - 2\sqrt{2}}\|q\|_{1,\rho}\right]\footnote{Here $[\cdot]$ denotes the integer part of a real number.}+1$$ and $$\tilde{\alpha}_n = \frac{3\sqrt{2}\pi}{n}\|q\|_{1,\rho}.$$
 The FD-method described by formulas \eqref{A_eq_1}, \eqref{A_eq_6}, \eqref{formula_for_u_j}, \eqref{Cauchy_kernel}, \eqref{TJ_eq_3} and \eqref{TJ_eq_4} converges to the eigensolution $(u_{n}(x); \lambda_{n})$ of problem \eqref{eq_1}, \eqref{eq_1_1} for all $n>n_{0}$. Furthermore, for the $n>n_{0}$ the following estimations of the method's convergence rate hold true:
 \begin{equation}\label{TJ_eq_25}
    \left\|u_{n}(x)-\stackrel{m}{u}_{n}\!\!(x)\right\|_{\infty, 1/\sqrt{\rho}}\leq \frac{\tilde{\alpha}_{n}^{m+1}}{(2m+3)\sqrt{\pi(m+2)}(1-\tilde{\alpha}_{n})},
 \end{equation}
 \begin{equation}\label{TJ_eq_26}
    \left|\lambda_{n}-\stackrel{m}{\lambda}_{n}\right|\leq \|q\|_{1,\rho}\frac{\tilde{\alpha}_{n}^{m}}{(2m+1)\sqrt{\pi (m+1)}(1-\tilde{\alpha}_{n})},
 \end{equation}
 where
 \begin{equation}\label{TJ_eq_27}
    \stackrel{m}{u}_{n}\!\!(x)=\sum\limits_{j=0}^{m} u_{n}^{(j)}(x),\quad \stackrel{m}{\lambda}_{n}=\sum\limits_{j=0}^{m} \lambda_{n}^{(j)}.
 \end{equation}
\end{theorem}

\section{The FD-method: software implementation}

In the section below we discuss the software implementation that was produced of the present method and describe explicitly the algorithm used in this implementation.

The software implementation was written in the Python programming language version 2.7 using the libraries NumPy, SciPy, mpmath and matplotlib. The use of the NumPy library has allowed us to have floating-point variables with up to quadruple precision\footnote{If the code called upon by SciPy and NumPy is compiled for the \texttt{x86\_64} architecture. For reasons to do the GCC compiler the same \texttt{numpy.longdouble} type we use results in 80-bit precision on 32-bit processors.}. We faced a technical problem when trying to compute the values of Legendre $Q_n$ function for an argument that's sufficiently close to $\pm1$ using SciPy's \texttt{lqmn} to circumvent which we had to resort to calling the corresponding function \texttt{legenq} of the mpmath library. This process involves converting the argument of \texttt{legenq} from the data type \texttt{numpy.longdouble} to \texttt{mpf} and back again with sufficient precision.

In the algorithm we use the {\it tanh rule} and {\it Stenger's formula} in order to approximate integration in (\ref{formula_for_lambda_j}), (\ref{formula_for_u_j}):

 \begin{equation}\label{int_ab}
\int_{a}^{b}f(x)dx=\int_{-\infty}^{+\infty}f\left(\frac{a+be^{t}}{1+e^{t}}\right)\frac{(b-a)dt}{(e^{-t/2}+e^{t/2})^{2}}\approx
 \end{equation}

$$
\approx h_{sinc}\sum_{i=-K}^{K}f\left(\frac{a+be^{ih_{sinc}}}{1+e^{ih_{sinc}}}\right)\frac{b-a}{(e^{-ih_{sinc}/2}+e^{ih_{sinc}/2})^{2}},
$$

 \begin{equation}\label{int_az}
\text{   }\int_{a}^{z_{j}}f(x)dx\approx h_{sinc}\sum_{i=-K}^{K}\delta_{j-i}^{(-1)}f\left(\frac{a+be^{ih_{sinc}}}{1+e^{ih_{sinc}}}\right)\frac{b-a}{(e^{-ih_{sinc}/2}+e^{ih_{sinc}/2})^{2}}
 \end{equation}
where
 $\delta_{i}^{(-1)}=\frac{1}{2}+\int_{0}^{i}\frac{\sin(\pi t)}{\pi t}dt,\  i=-2K\ldots 2K$,
$h_{sinc}=\sqrt{\frac{2\pi}{K}}$.

Below we also use  the following auxiliary notation:

\begin{equation}
z_{i}=\frac{a+be^{h_{sinc}i}}{1+e^{h_{sinc}i}},\ \mu_{i}=\frac{b-a}{(e^{-ih_{sinc}/2}+e^{ih_{sinc}/2})^{2}},
\end{equation}
and refer to $A^{-1}$ as defined in (\ref{A_inv}).

In order to measure how close an obtained approximation is to the exact solution we used the functional

$$
\text{   }\overset{m}{\eta}_{n}=\left[\int_{-1}^{1}\left[(1-x^{2})\frac{d \overset{m}u_{n}(x)}{dx}+\int_{-1}^{x}\left(\overset{m}\lambda_{n}-q(\xi)\right)\overset{m}{u}_{n}(\xi)d\xi\right]^{2}dx\right]^{\frac{1}{2}}
$$
referred to in the algorithm as the \textit{residual.}

The developed software library implements the capacity to subdivide the interval $(a, b)$ on which numerical integration takes place into subintervals $(a = x_0, x_1), \dots, (x_{N-1}, x_{N} = b)$ in a uniform as well as a non-uniform manner. A separate set of $z_i, \mu_i$ is generated for each $(x_{i-1}, x_{i}), i \in \{0, 1, \dots, N\}$ in that case. Since $q(x)$ is sampled at the points $z_i$, which are at their densest at the ends of the interval, one could benefit from subdividing the interval at the singularity points of $q(x)$. For the sake of simplicity we shall omit this detail in the description of the algorithm that follows.

As the values of $\delta_{i}^{(-1)}$ do not depend on $q(x)$ or how the interval is subdivided they were precomputed and stored in a file to be loaded by the library at runtime.

Note: when in the algorithm we say ``$F[i][j]$'' we refer to a particular element of the two-dimensional array $F$ that has the index $i, j$. However, when we refer to ``$F[i]$'' what we mean is the values $F[i][-K], F[i][-K + 1], \dots, F[i][K]$ taken as a one-dimensional array.

\newcommand{\IntAB}{\texttt{IntAB}}
\newcommand{\IntAZ}{\texttt{IntAZ}}

The main computing routine is described in Algorithm \ref{alg_main}. It references the subroutines $\IntAB$ and $\IntAZ$ defined in Algorithms \ref{alg_intab}, \ref{alg_intaz}.

\begin{algorithm}
    \SetAlgoLined
    \caption{IntAB(values)}
	\label{alg_intab}
    \KwData{$values$, $h_{sinc}$, $z_j$, $\mu_j$}
    \KwResult{$s$}
    \Begin{
    	$s := 0$\;
	     \For{$j := -K \ldots K$}{
	     	$p := \mu_j$\;
		    \ForEach{v \textbf{in} values}{
		          	\uIf{$v$ is a function}{
		          		$p := p \, v(z_j)$\;
		          	}
		          	\Else{\tcp{v is an array}
		          		$p := p \,  v[j]$\;		          	
		          	}
		    }
		    $s := s + p$\;
	    }
	  $s := h_{sinc} \, s$\;
    }
\end{algorithm}

\begin{algorithm}
    \SetAlgoLined
    \caption{IntAZ(j;values)}
	\label{alg_intaz}
    \KwData{$j$, $values$, $h_{sinc}$, $z_j$, $\mu_i$, $\delta^{(-1)}_i$}
    \KwResult{$r$}
    \Begin{
    	$s := 0$\;
	     \For{$i := -K \ldots K$}{
	     	$p := \mu_i \delta^{(-1)}_{j - i}$\;
		    \ForEach{v \textbf{in} values}{
		          	\uIf{$v$ is a function}{
		          		$p := p \, v(z_i)$\;
		          	}
		          	\Else{\tcp{v is an array}
		          		$p := p \,  v[i]$\;		          	
		          	}
		    }
		    $s := s + p$\;
	    }
	  $s := h_{sinc} \, s$\;
    }
\end{algorithm}

\begin{algorithm}
    \SetAlgoLined
    \caption{Main}
	\label{alg_main}
    \KwData{$n$ --- the number of the eigenvalue we want to find, $m$ --- the order of the FD-method (the number of steps taken), $K$,  $h_{sinc}$, $z_i$, $\mu_i$, $\delta^{(-1)}_i$}
    \KwResult{$\overset{m}{\lambda_n}$, $\overset{m}{\eta_n}$, $\overset{m}{u_n}(x)$, $\frac{d\overset{m}{u}_n}{dx}(x)$, $\left\{ \left\Vert u_{n}^{(i)}(x)\right\Vert \right\} _{i=0}^{m}$}
    \Begin{
		\tcp{We initialize L as a one-dimensional array of 2K + 1 zeros and F, U and DU as two-dimensional arrays of 2K + 1 by 2K + 1 zeros.}
        $L := zeros(-K \ldots K)$\;	
        $F, U, DU := zeros(-K \ldots K, -K \ldots K)$\;
        $L[0] = n (n + 1)$\;
	    \For{$i := -K \ldots K$}{
        	$U[0][i] = P_n(x)$\;
        	$DU[0][i] = dP_n(x)$\;         	
	    }
        \For{$d := 1,2 \ldots m$}{
            \tcp{Compute the correction for the eigenvalue}
        	$L[d] := A^{-2} \IntAB(U[0], U[d - 1], q)$\;
            \tcp{Compute F}
	        \For{$i := -K \ldots K$}{
		        $F[d][i] := U[d - 1][i]\, q(z_i)$\;	
		        \For{$j := 0  \ldots d - 1$}{
					$F[d][i] := F[d][i] - L[d - j]\,U[j][i]$\;
		        }
		    }                  	
            \tcp{Compute the correction for the eigenfunction}
	        \For{$i := -K \ldots K$}{
            	$U[d][i] := Q_n(z_i) \IntAZ(i; F[d], P_n) - P_n(z_i) \IntAZ(i; F[d], Q_n)$\;
            	$DU[d][i] := dQ_n(z_i) \IntAZ(i; F[d], P_n) - dP_n(z_i)   \IntAZ(i; F[d], Q_n)$\;            	
		    }
            \tcp{Orthogonality}
            $I = A^{-2} \IntAB (U[d], U[0])$\;
	        \For{$i := -K \ldots K$}{
            	$U[d][i] := U[d][i] - I\,U[0][i]$\;
            	$DU[d][i] := DU[d][i] - I\, DU[0][i]$\;            	
		    }
            \tcp{Compute the residual}
            $\texttt{CompRes}$\;
        }
        $\overset{m}{\lambda_n} := \sum_{i=0}^{m} L[i]$;
    }
\end{algorithm}

\section{Numerical experiments}

Using the above algorithm we applied the FD-method to problem (\ref{eq_1}), (\ref{eq_1_1}) with the potential

$$
q(x)=\ln\left(\left|\left(\frac{5}{12}-x\right)\left(\frac{1}{3}+x\right)\right|\right).
$$

First, the software was run to approximate the value of $\lambda_0$ with $m = 60$ steps of the FD-method to demonstrate the rate of convergence. In this and subsequent runs the quadrature formulas (\ref{int_ab}), (\ref{int_az}) had $K$ at $250$ and for numerical integration $(-1, 1)$ was subdivided into four subintervals using the set of points $\{-1, -\frac{1}{3}, 0, \frac{5}{12}, 1\}$. The importance of using the latter kind of subdivision is illustrated below.

\begin{figure}[h]
	\centering
	\begin{minipage}{0.4\textwidth}
		\centering
		 \includegraphics[width=0.99\linewidth]{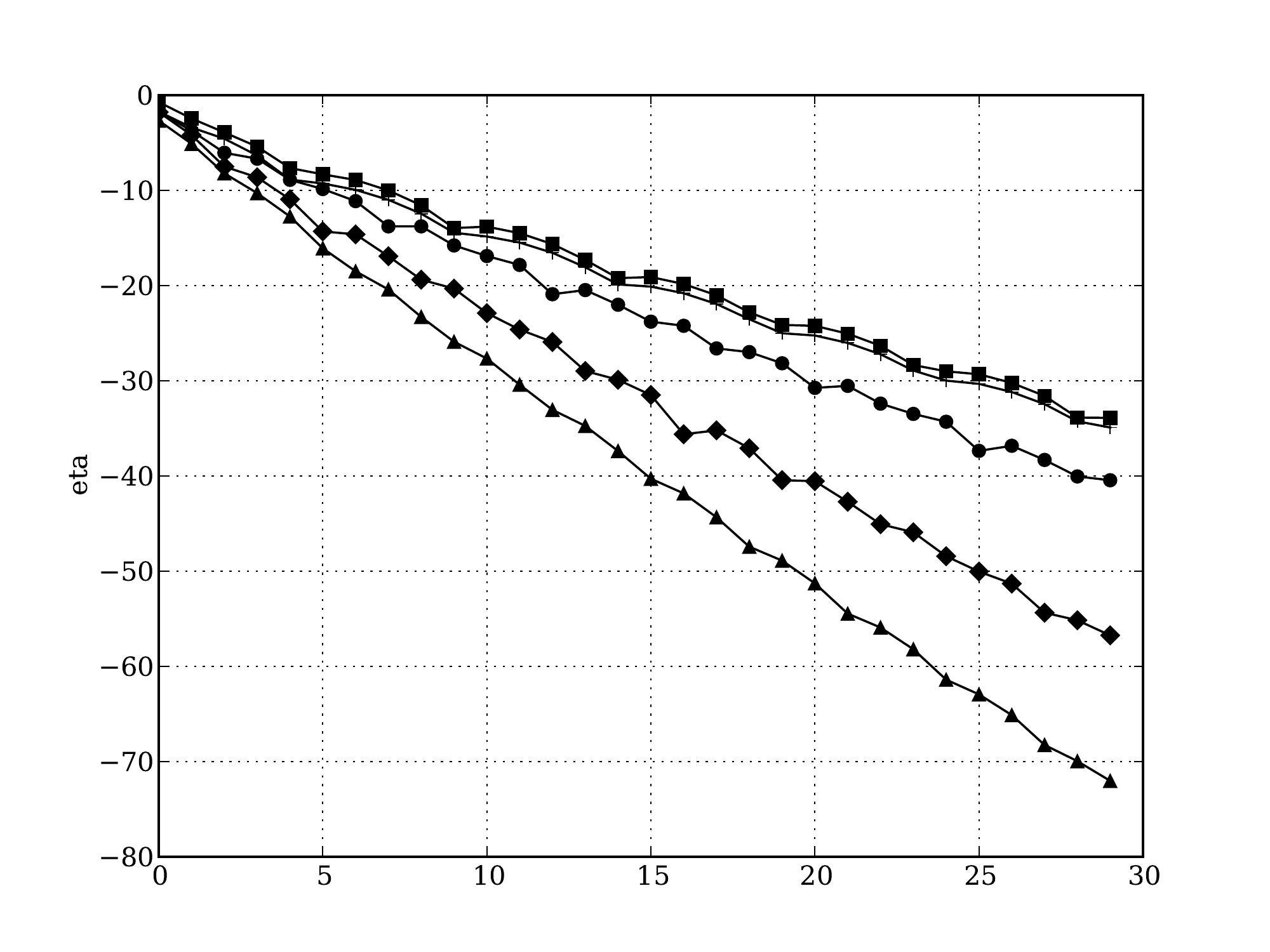}
		$\ln\left(\overset{m}{\eta}_{n}\right)$	
	\end{minipage}
	\begin{minipage}{0.1\textwidth}
		\centering
    	\includegraphics[width=0.6\linewidth]{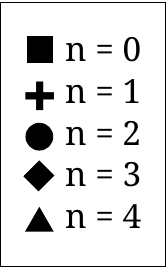}	
	\end{minipage}	
	\begin{minipage}{0.4\textwidth}
		\centering
		 \includegraphics[width=0.99\linewidth]{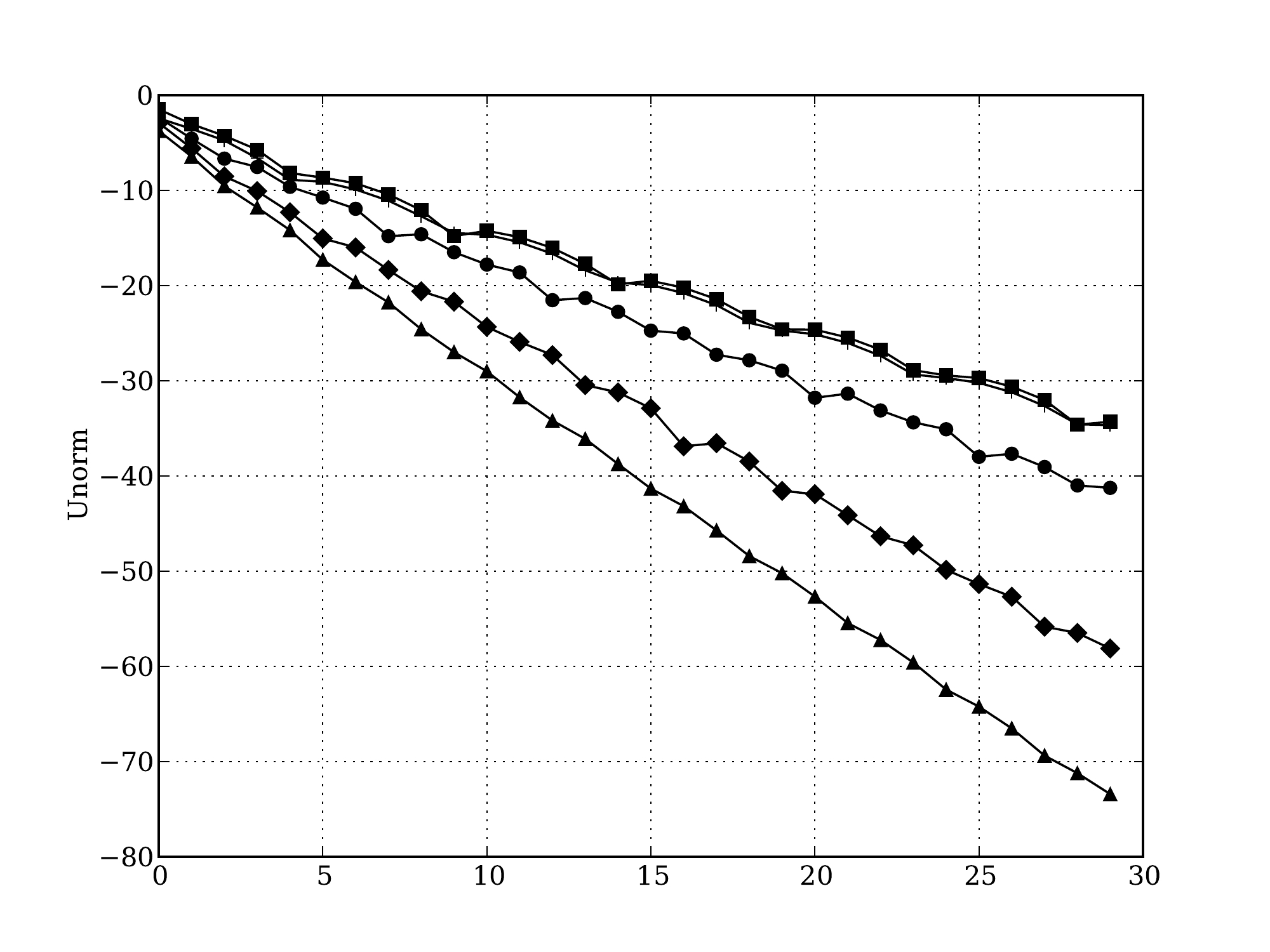}
		$\ln\left(\left\Vert u_{n}^{(i)}(x)\right\Vert\right)$
	\end{minipage}
	\caption{A log-scale graph that shows the convergence rate for $\lambda_0, \dots, \lambda_4$.
	}
	\label{graphs}	
\end{figure}

Figure \ref{graphs} illustrates how the convergence rate of the FD-method increases exponentially for each subsequent eigenvalue $\lambda_n$.

Solving the same problem was attempted using the well-known SLEIGN2 software package. The rightmost columns of Tables \ref{conv_lambda0}, \ref{results_n05_table} show the margin of error in the results it produces compared to the present implementation of the FD-method.

\begin{table}[h]
\caption{The results obtained using SLEIGN2}
\begin{tabular}{|l|l|l|l|}
\hline
{$n$}&	{$\lambda_{n,sl2}$}&	{$TOL$}&	{{\it IFLAG}}	\\
\hline
{$0$}&	{$-1.98326983D+00$}&	{$0.46748D-08$}&	{$1$}	\\
\hline
{$1$}&	{$ 0.855187683D+00$}&	{$0.73426D-07$}&	{$1$}	\\
\hline
{$2$}&	{$ 0.489606686D+01$}&	{$0.35447D-07$}&	{$1$}	\\
\hline
{$3$}&	{$ 0.104183770D+02$}&	{$0.40228D-07$}&	{$1$}	\\
\hline
{$4$}&	{$ 0.188163965D+02$}&	{$0.61329D-11$}&	{$1$}	\\
\hline
\end{tabular}
\end{table}

\begin{table}[h]
\caption{Convergence for the eigenvalue $\lambda_0$}
\label{conv_lambda0}
\begin{tabular}{|l|l|l|l|l|l|}

\hline
$m$ & $\overset{m}{\lambda_0}$ & $\left\Vert u_{0}^{(m)}(x)\right\Vert$ &  $\overset{m}{\eta}_{0}$ & $|\overset{m}{\lambda} _0-\lambda_{0,sl2}|$\\
\hline
0&-1.8538570587&0.2270941786&0.4851738751&0.1294127713\\
1&-2.0002817053&0.0478946893&0.0863600316&0.0170118753\\
2&-1.9826820263&0.0140616365&0.0200439342&0.0005878037\\
3&-1.9827492251&0.0032752573&0.0044828399&0.0005206049\\
4&-1.9832100727&0.000281665&0.0004743141&0.0000597573\\
5&-1.9831500665&0.0001734894&0.0002452252&0.0001197635\\
6&-1.9831433619&9.61416137299e-05&0.0001358145&0.0001264681\\
7&-1.9831424182&2.99030462249e-05&4.5191830517e-05&0.0001274118\\
8&-1.9831451284&5.71179849936e-06&9.49092804135e-06&0.0001247016\\
9&-1.9831441732&3.7195952769e-07&8.68014240854e-07&0.0001256568\\
\hline
\end{tabular}
\end{table}

\begin{table}[h]
\begin{tabular}{|l|l|l|l|l|l|}
\hline
$m$ & $\overset{m}{\lambda_0}$ & $\left\Vert u_{0}^{(m)}(x)\right\Vert$ &  $\overset{m}{\eta}_{0}$ & $|\overset{m}{\lambda} _0-\lambda_{0,sl2}|$\\
\hline
50&-1.983144271&3.6458910063e-24&5.42202004605e-24&0.000125559\\
51&-1.983144271&1.45758164365e-24&2.17584093998e-24&0.000125559\\
52&-1.983144271&3.59257473326e-25&5.46661601878e-25&0.000125559\\
53&-1.983144271&2.29601032831e-26&5.40330904597e-26&0.000125559\\
54&-1.983144271&4.24716018236e-26&6.37606340182e-26&0.000125559\\
55&-1.983144271&2.6000804557e-26&3.86695510098e-26&0.000125559\\
56&-1.983144271&9.47336776012e-27&1.41695598924e-26&0.000125559\\
57&-1.983144271&2.02349833941e-27&3.116321958e-27&0.000125559\\
58&-1.983144271&1.81719782215e-28&3.7965999043e-28&0.000125559\\
59&-1.983144271&3.43816989604e-28&5.133026339e-28&0.000125559\\
60&-1.983144271&1.8365375822e-28&2.7327040704e-28&0.000125559\\
\hline
\end{tabular}
\end{table}

Computations for further eigenvalues were also performed and compared to the results from SLEIGN2 (see Tables \ref{results_n04_table}, \ref{results_n05_table}).

\begin{table}[h]
\caption{The values obtained for $\lambda_0,\dots,\lambda_4$ at $m = 30$}
\label{results_n04_table}
\begin{tabular}{|l|l|l|l|l|l|}
\hline
$n$ & $\overset{m}{\lambda_n}$ & $|\lambda^{(m)}_n|$ & $\left\Vert u_{n}^{(m)}(x)\right\Vert$ \\
\hline
0 & -1.98314427097744064 & 1.46303698262e-17 & 1.26598694672e-15
\\
1 & 0.857270328373118208 & 1.63565545758e-17 & 8.83118381572e-16
\\
2 & 4.893950682679907660 & 1.72618520779e-18 & 1.22013435336e-18
\\
3 & 10.42051129625743390 & 5.71577711655e-26 & 4.58227541331e-25
\\
4 & 18.81639652150898795 & 1.30790575077e-32 & 5.43628701044e-32
\\
\hline
\end{tabular}
\end{table}

\begin{table}[h]
\caption{Accuracy results for $\lambda_0,\dots,\lambda_4$ at $m = 30$}
\label{results_n05_table}
\begin{tabular}{|l|l|l|l|l|l|}
\hline
$n$ & $\overset{m}{\lambda_n}$  &  $\overset{m}{\eta}_{n}$ & $|\overset{m}{\lambda} _0-\lambda_{0,sl2}|$\\
\hline
0 & -1.98314427097744064  & 1.9052706379e-15 & 0.000125559
\\
1 & 0.857270328373118208  & 6.92114145514e-16 & 0.0020826454
\\
2 & 4.893950682679907660  & 2.72086325283e-18 & 0.0021161773
\\
3 & 10.42051129625743390  & 2.28096722974e-25
 & 0.0021342963
\\
4 & 18.81639652150898795  & 5.26360265358e-32 & 0.0000000215\\
\hline
\end{tabular}
\end{table}

For the eigenvalue $\lambda_0$ and $m = 30$ the choice of subdivision mattered significantly. Numerical experiments show that the subdivision points are best placed near the singularities of $q(x)$ (see Table \ref{results_n06_table}).

\begin{table}[h]
\caption{The values obtained for $\lambda_0$ at $m = 30$ with different subdivisions.}
\label{results_n06_table}\begin{center}
\begin{tabular}{|l|l|l|l|l|l|}
\hline
Subdivision & $\overset{m}{\lambda_0}$ & $|\overset{m}{\lambda} _0-\lambda_{0,sl2}|$\\
\hline
$N = 1$, none & -1.9318815213501200317 & 0.051388309 \\
\hline
$N = 4$, uniform & -1.9776298960768203497 & 0.005639934 \\
\hline
$N = 4$, $\{-1, -\frac{1}{3}, 0, \frac{5}{12}, 1\}$ & -1.9831442710817836887 & 0.000125559\\
\hline
\end{tabular}
\end{center}
\end{table}

\section{Conclusions}

The article lays out the structure of and provides a theoretical justification for the FD-method as applied to solving the Sturm-Liouville problem (\ref{eq_1}), (\ref{eq_1_1}). In Theorem \ref{conv_theorem} convergence is proven for the case when $q(x)$ satisfies condition  (\ref{eq_1_2}) and estimates for the convergence rate are given explicitly.

Special attention should also be drawn to Theorem \ref{Elegant_theorem}. The authors were unable to find analogous results in the existing literature. To their best knowledge the theorem  and its proof constitute a novel and original result.

The presented method suggests at least two ways for further refinement. First, by considering a separate approximation  for the potential on each subinterval of $[-1, 1]$ (as done for piecewise continuous potential problems in \cite{Dragunov_Math_Comp}). Second, by modifying the algorithm for concurrent computation. The authors hope to explore these possibilities in future publications.

The algorithm was implemented in software as a function library (a Python module). The implementation can be integrated into larger systems or used as is in applied sciences. The source code for the function library  along with example Python code that uses it can be obtained from \url{https://github.com/imathsoft/legendrefdnum}.

\bibliographystyle{plain}

\bibliography{references_stat}

\end{document}